\documentstyle[thmsa,a4,12pt,sw20lart]{article}

\input{tcilatex}
\begin{document}

\author{Teodor Oprea \\
University of Bucharest,\\
Faculty of Mathematics and Informatics\\
$14$ Academiei St., code 010014, Bucharest, Romania\\
e-mail: teodoroprea@yahoo.com}
\title{Chen's inequality in Lagrangian case}
\date{}
\maketitle

\begin{center}
\medskip {\bf Abstract}
\end{center}

\medskip

In the theory of submanifolds, the following problem is fundamental: {\sl to
establish simple relationships between the main intrinsic invariants and the
main extrinsic invariants of the submanifolds. }The basic relationships
discovered until now [1{\bf ,} 2, 3, 4] are inequalities. To analyze these
problems, we follow the idea of C. Udri\c{s}te [6{\bf ,} 7, 8{\bf ] }that
the method of constrained extremum is a natural way to prove geometric
inequalities. We improve the Chen's inequality which characterizes a totally
real submanifold of a complex space form. For that we suppose that the
submanifold is Lagrangian and we formulate and analyze a suitable
constrained extremum problem.

\medskip

2000{\bf \ }{\it Mathematics subject classification: }53C21, 53C24, 53C25,
49K35.


\bigskip

\begin{center}
{\bf 1.Optimizations on Riemannian submanifolds}\ 
\end{center}

\medskip

Let $(N,\widetilde{g}$ $)$ be a Riemannian manifold of dimension $m,$ $M$ be
a Riemannian submanifold of it, $g$ be the metric induced on $M$ by $%
\widetilde{g}$ and $f:N\rightarrow $ $R$ be a differentiable function.

In [5] we considered the constrained extremum problem

\begin{center}
(1) $\stackunder{x\in M}{\min }f(x).$
\end{center}

The first result is

\ 

T{\footnotesize HEOREM} 1.{\sl \ If }$x_{0}\in M${\sl \ is the solution of
the problem }$(1)${\sl , then}

i){\sl \ }$($grad f$)(x_{0})${\sl \ }$\in T_{x_{0}}^{\perp }M,$

ii){\sl \ the bilinear form}

\begin{center}
$\alpha {\sl \ :\ }T_{x_{0}}M\times T_{x_{0}}M{\sl \ }\rightarrow {\sl \ }R,$

$\alpha (X,Y)=$Hess$_{f}(X,Y)+\widetilde{g\text{ }}(h(X,Y),($grad f$%
)(x_{0})),${\sl \ }
\end{center}

{\sl is positive semidefinite, where }$h${\sl \ is the second fundamental
form of the submanifold }$M${\sl \ in }$N.$

\ 

We shall use this theorem in order to find a inequality verified by the
Chen's invariant of a Lagrangian submanifold in a complex space form.

\begin{center}
\bigskip \bigskip

{\bf 2. Estimation of Chen's invariant of a Lagrangian submanifold in a
complex space form}
\end{center}

\medskip

Let $(M,g)$ be a Riemannian manifold of dimension $n$, and $x$ a point in $%
M. $ We consider the orthonormal frame $\{e_{1},e_{2},...,e_{n}\}$ in $%
T_{x}M.$ The scalar curvature at $x$ is defined by

\[
\tau =\dsum\limits_{1\leq i<j\leq n}R(e_{i},e_{j},e_{i},e_{j}), 
\]
where $R$ is the Riemann curvature tensor of $(M,g).$

We denote $\delta _{M}=\tau -\min (k)$, where $k$ is the sectional curvature
at the point $x.$ The invariant $\delta _{M}$ is called the Chen's invariant
of Riemannian manifold $(M,g).$

Let $(\widetilde{M},\widetilde{g},J)$ be a K\"{a}hler manifold of real
dimension $2m.$ A submanifold $M$ of dimension $n$ of $(\widetilde{M},%
\widetilde{g},J)$ is called a totally real submanifold if for any point $x$
in $M$ the relation $J(T_{x}M)\subset T_{x}^{\perp }M$ holds.

If, in addition, $n=m,$ then $M$ is called Lagrangian submanifold. For a
Lagrangian submanifold, the relation $J(T_{x}M)=T_{x}^{\perp }M$ occurs.

A K\"{a}hler manifold with constant holomorphic sectional curvature is
called a complex space form and is denoted by $\widetilde{M}(c)$. The
Riemann curvature tensor $\widetilde{R}$ of $\widetilde{M}(c)$ satisfies the
relation

$\widetilde{R}(X,Y)Z=\frac{c}{4}\{\widetilde{g}(Y,Z)X-\widetilde{g}(X,Z)Y+%
\widetilde{g}(JY,Z)JX-\widetilde{g}(JX,Z)JY+2\widetilde{g}(X,JY)JZ\}.$

A totally real submanifold of real dimension $n$ in a complex space form $%
\widetilde{M}(c)$ of real dimension $2m$ verifies a Chen's inequality

\[
\delta _{M}\leq \frac{n-2}{2}\{\frac{n^{2}}{n-1}\left\| H\right\| ^{2}+(n+1)%
\frac{c}{4}\}, 
\]
where $H$ is the mean curvature vector of the Riemannian submanifold $M$ of $%
\widetilde{M}(c).$

R{\footnotesize EMARK}{\bf . }i) If $M$ is a totally real submanifold of
real dimension $n$ in a complex space form $\widetilde{M}(c)$ of real
dimension $2m,$ then

\begin{center}
\medskip $A_{JY}X=-Jh(X,Y)=A_{JX}Y,$ $\forall $ $X,Y\in {\cal X}(M).$
\end{center}

ii) Let $m=n$ ($M$ is Lagrangian in $\widetilde{M}(c)$). If we consider the
point $x\in M$, the orthonormal frames $\{e_{1},...,e_{n}\}$ in $T_{x}M$ and 
$\{Je_{1},...,Je_{n}\}$ in $T_{x}^{\perp }M$, then

\medskip 
\[
h_{jk}^{i}=h_{ik}^{j},\forall i,j,k\in \overline{1,n}, 
\]
where $h_{jk}^{i}$ is the component after $Je_{i}$ of the vector $%
h(e_{j},e_{k}).$

\ 

With these ingredients we prove the next result which can be regarded as an
obstruction to Lagrangian isometric immersions of a Riemannian manifold into
a complex space form.

\ 

T{\footnotesize HEOREM 2.}{\bf \ }{\sl Let }$M${\sl \ be a Lagrangian
submanifold in complex space form }$\widetilde{M}(c)${\sl \ of real
dimension }$2n${\sl , }$n\geq 3.${\sl \ Then }

\begin{center}
$\delta _{M}\leq \frac{(n-2)(n+1)}{2}\frac{c}{4}+\frac{n^{2}}{2}\frac{2n-3}{%
2n+3}\left\| H\right\| ^{2}.$
\end{center}

\ 

P{\footnotesize ROOF}{\bf . }We consider the point{\bf \ }$x\in M$, the
orthonormal frames $\{e_{1},...,e_{n}\}$ in $T_{x}M$ and $%
\{Je_{1},...,Je_{n}\}$ in $T_{x}^{\perp }M$, $\{e_{1},e_{2}\}$ being an
orthonormal frame in the $2-$ plane which minimize the sectional curvature
at the point $x$ in $T_{x}M.$

By using Gauss's equation $\widetilde{R}(X,Y,Z,U)=$ $\widetilde{R}(X,Y,Z,U)-$%
\newline
$-g(h(X,Z),h(Y,U))+g(h(X,U),h(Y,Z)),\forall $ $X,Y,Z,U\in {\cal X}(M)$ and
the fact that $\widetilde{M}(c)$ is a complex space form, we obtain\newline
(1) $\tau =\frac{n(n-1)}{2}\frac{c}{4}+\dsum\limits_{r=1}^{n}\dsum\limits_{1%
\leq i<j\leq
n}h_{ii}^{r}h_{jj}^{r}-\dsum\limits_{r=1}^{n}\dsum\limits_{1\leq i<j\leq
n}(h_{ij}^{r})^{2}$,\newline
(2) $R(e_{1},e_{2},e_{1},e_{2})=\frac{c}{4}+\dsum%
\limits_{r=1}^{n}h_{11}^{r}h_{22}^{r}-\dsum%
\limits_{r=1}^{n}(h_{12}^{r})^{2}. $

By subtracting the relations (1) and (2), we find\newline
(3) $\delta _{M}=\frac{(n-2)(n+1)}{2}\frac{c}{4}+\dsum\limits_{r=1}^{n}(%
\dsum\limits_{3\leq j\leq
n}(h_{11}^{r}+h_{22}^{r})h_{jj}^{r}+\dsum\limits_{3\leq i<j\leq
n}h_{ii}^{r}h_{jj}^{r}-$\newline
$-\dsum\limits_{3\leq j\leq n}(h_{1j}^{r})^{2}-\dsum\limits_{2\leq i<j\leq
n}(h_{ij}^{r})^{2}).$

By using the symmetry in the three indexes of $h_{ij}^{k},$ we can write%
\newline
(4) $\delta _{M}\leq \frac{(n-2)(n+1)}{2}\frac{c}{4}+\dsum\limits_{r=1}^{n}(%
\dsum\limits_{3\leq j\leq
n}(h_{11}^{r}+h_{22}^{r})h_{jj}^{r}+\dsum\limits_{3\leq i<j\leq
n}h_{ii}^{r}h_{jj}^{r})-\newline
-\dsum\limits_{3\leq j\leq n}(h_{1j}^{1})^{2}-\dsum\limits_{3\leq j\leq
n}(h_{1j}^{j})^{2}-\dsum\limits_{2\leq i<j\leq
n}(h_{ij}^{i})^{2}-\dsum\limits_{2\leq i<j\leq n}(h_{ij}^{j})^{2}=$\newline
$=\frac{(n-2)(n+1)}{2}\frac{c}{4}+\dsum\limits_{r=1}^{n}(\dsum\limits_{3\leq
j\leq n}(h_{11}^{r}+h_{22}^{r})h_{jj}^{r}+\dsum\limits_{3\leq i<j\leq
n}h_{ii}^{r}h_{jj}^{r})-$\newline
$-\dsum\limits_{3\leq j\leq n}(h_{11}^{j})^{2}-\dsum\limits_{3\leq j\leq
n}(h_{jj}^{1})^{2}-\dsum\limits_{2\leq i<j\leq
n}(h_{ii}^{j})^{2}-\dsum\limits_{2\leq i<j\leq n}(h_{jj}^{i})^{2}=$\newline
$=\frac{(n-2)(n+1)}{2}\frac{c}{4}+\dsum\limits_{r=1}^{n}(\dsum\limits_{3\leq
j\leq n}(h_{11}^{r}+h_{22}^{r})h_{jj}^{r}+\dsum\limits_{3\leq i<j\leq
n}h_{ii}^{r}h_{jj}^{r})-$\newline
$-\dsum\limits_{3\leq j\leq n}(h_{11}^{j})^{2}-\dsum\limits_{3\leq j\leq
n}(h_{jj}^{1})^{2}-\dsum\limits\Sb i,j\in \overline{2,n}  \\ i\neq j  \endSb %
(h_{jj}^{i})^{2}.$

\ 

Let us consider the quadratic forms $f_{1},$ $f_{2},$ $f_{r}:R^{n}%
\rightarrow R,$ $r\in \overline{3,n}$ defined respectively by

$\ $

$f_1(h_{11}^1,h_{22}^1,...,h_{nn}^1)=\dsum\limits_{3\leq j\leq
n}(h_{11}^1+h_{22}^1)h_{jj}^1+\dsum\limits_{3\leq i<j\leq
n}h_{ii}^1h_{jj}^1-\dsum\limits_{3\leq j\leq n}(h_{jj}^1)^2,$

$f_2(h_{11}^2,h_{22}^2,...,h_{nn}^2)=\dsum\limits_{3\leq j\leq
n}(h_{11}^2+h_{22}^2)h_{jj}^2+\dsum\limits_{3\leq i<j\leq
n}h_{ii}^2h_{jj}^2-\dsum\limits_{3\leq j\leq n}(h_{jj}^2)^2,$

$f_r(h_{11}^r,h_{22}^r,...,h_{nn}^r)=\dsum\limits_{3\leq j\leq
n}(h_{11}^r+h_{22}^r)h_{jj}^r+\dsum\limits_{3\leq i<j\leq
n}h_{ii}^rh_{jj}^r- $

\begin{center}
$-(h_{11}^{r})^{2}-\dsum\limits\Sb j\in \overline{2,n}  \\ j\neq r  \endSb %
(h_{jj}^{r})^{2}.$
\end{center}

We need the maximum of $f_{1}$ and $f_{3}$. For $f_{2}$ and $f_{r}$, $r\in 
\overline{4,n}$, we can solve similar problems.

We start with the problem 
\[
\max f_{1}, 
\]

\[
\text{subject to }P:h_{11}^{1}+h_{22}^{1}+...+h_{nn}^{1}=k^{1}, 
\]
where $k^{1}$ is a real constant.

The first three partial derivatives of the function $f_{1}$ are\newline
(5) $\frac{\partial f_{1}}{\partial h_{11}^{1}}=\dsum\limits_{3\leq j\leq
n}h_{jj}^{1},$\newline
(6) $\frac{\partial f_{1}}{\partial h_{22}^{1}}=\dsum\limits_{3\leq j\leq
n}h_{jj}^{1},$\newline
(7) $\frac{\partial f_{1}}{\partial h_{33}^{1}}=h_{11}^{1}+h_{22}^{1}+\dsum%
\limits_{4\leq j\leq n}h_{jj}^{1}-2h_{33}^{1}.$

As for a solution $(h_{11}^{1},h_{22}^{1},...,h_{nn}^{1})$ of the problem in
question, the vector\newline
$($grad) $(f_{1})$ is normal at $P$, from (5), (6) and (7) we obtain\newline
(8) $h_{11}^{1}+h_{22}^{1}=3h_{jj}^{1}=3a^{1},$ $\forall $ $j\in \overline{%
3,n}.$

By using the relation $%
h_{11}^{1}+h_{22}^{1}+h_{33}^{1}+...+h_{nn}^{1}=k^{1}, $ from (8) we obtain $%
3a^{1}+(n-2)a^{1}=k^{1}.$ Consequently\newline
(9) $a^{1}=\frac{k^{1}}{n+1}.$\ 

As $f_{1\text{ }}$is obtained from the function studied in Chen's inequality
(see [5]) by subtracting some square terms, $f_{1}\left| P\right. $ will
have the Hessian seminegative definite. Consequently the point $%
(h_{11}^{1},h_{22}^{1},...,h_{nn}^{1})$ given by the relations (8) and (9)
is a global maximum point, and hence\newline
(10) $f_{1}\leq 3a^{1}(n-2)a^{1}+C_{n-2}^{2}(a^{1})^{2}-(n-2)(a^{1})^{2}=%
\frac{(a^{1})^{2}}{2}(n+1)(n-2).$

\ From (9) and (10), it follows\newline
(11) $f_{1}\leq \frac{(k^{1})^{2}}{2}\frac{n-2}{n+1}=\frac{(n)^{2}}{2}\frac{%
n-2}{n+1}(H^{1})^{2}.$

Similarly we obtain\newline
(12) $f_{2}\leq \frac{(n)^{2}}{2}\frac{n-2}{n+1}(H^{2})^{2}.$

Further on, we shall consider the problem

\[
\max f_3, 
\]
\[
\text{subject to }P:h_{11}^3+h_{22}^3+...+h_{nn}^3=k^3, 
\]
where $k^3$ is a real constant.

The first four partial derivatives of the function $f_{3}$ are\newline
(13) $\frac{\partial f_{3}}{\partial h_{11}^{3}}=\dsum\limits_{3\leq j\leq
n}h_{jj}^{3}-2h_{11}^{3},$\newline
(14) $\frac{\partial f_{3}}{\partial h_{22}^{3}}=\dsum\limits_{3\leq j\leq
n}h_{jj}^{3}-2h_{22}^{3},$\newline
(15) $\frac{\partial f_{3}}{\partial h_{33}^{3}}=h_{11}^{3}+h_{22}^{3}+\dsum%
\limits_{4\leq j\leq n}h_{jj}^{3},$\newline
(16) $\frac{\partial f_{3}}{\partial h_{44}^{3}}=h_{11}^{3}+h_{22}^{3}+\dsum%
\limits\Sb 3\leq j\leq n  \\ j\neq 4  \endSb h_{jj}^{3}-$\ $2h_{44}^{3}.$

For a solution $(h_{11}^{1},h_{22}^{1},...,h_{nn}^{1})$ of the problem in
question, the vector $($grad$)(f_{3})$ is colinear to $(1,1,...,1)$.

By using (13), (14), (15) and (16) we obtain\newline
(17) $h_{11}^{3}=h_{22}^{3}=3a^{3},$\newline
(18) $h_{33}^{3}=12a^{3},$\newline
(19) $h_{jj}^{3}=4a^{3},$ $\forall $ $j\in \overline{4,n}.$

As $h_{11}^{3}+h_{22}^{3}+h_{33}^{3}+...+h_{nn}^{3}=k^{3},$ from (17), (18)
and (19), one gets\newline
(20) $a^{3}=\frac{k^{3}}{4n+6}.$

With an argument similar to those in the previous problem we obtain that the
point $(h_{11}^{3},h_{22}^{3},...,h_{nn}^{3})$ given by the relations (17),
(18), (19) and (20) is a global maximum point. Therefore\newline
(21) $f_{3}\leq
6a^{3}12a^{3}+6a^{3}(n-3)4a^{3}+12b(n-3)4a^{3}+C_{n-3}^{2}16(a^{3})^{2}-$

\begin{center}
$-18(a^{3})^{2}-(n-3)16(a^{3})^{2}=2(a^{3})^{2}(2n-3)(2n+3).$
\end{center}

From (20) and (21) we obtain $f_{3}\leq \frac{(k^{3})^{2}}{2}\frac{2n-3}{2n+3%
}=\frac{n^{2}}{2}\frac{2n-3}{2n+3}(H^{3})^{2}.$\ 

Similarly one gets\newline
(22) $f_{r}\leq \frac{n^{2}}{2}\frac{2n-3}{2n+3}(H^{r})^{2},$ $\forall $ $%
r\in \overline{3,n}.$

As $\frac{n-2}{n+1}<\frac{2n-3}{2n+3}$, $\forall $ $n\geq 3$, from (11),
(12) and (22), it follows\newline
(23) $f_{r}\leq \frac{n^{2}}{2}\frac{2n-3}{2n+3}(H^{r})^{2},$ $\forall $ $%
r\in \overline{1,n}.$

By using the relations (4) and (23) we have\newline
(24) $\delta _{M}\leq \frac{(n-2)(n+1)}{2}\frac{c}{4}+\frac{n^{2}}{2}\frac{%
2n-3}{2n+3}\dsum\limits_{r=1}^{n}(H^{r})^{2}=\frac{(n-2)(n+1)}{2}\frac{c}{4}+%
\frac{n^{2}}{2}\frac{2n-3}{2n+3}\left\| H\right\| ^{2}.$

\ 

\medskip

R{\footnotesize EMARK}{\bf . }B.Y. Chen, F. Dillen, L. Verstraelen, L.
Vrancken showed in [4] that an Lagrangian submanifold, of real dimension $2n$%
, $n\geq 3,$ of a complex space form $\widetilde{M}(c)$, satisfying the
equality

\medskip 
\[
(25)\text{ }\delta _{M}=\frac{n-2}{2}\{\frac{n^{2}}{n-1}\left\| H\right\|
^{2}+(n+1)\frac{c}{4}\}, 
\]
is minimal.

Now, this result is an immediate consequence of the previous inequality.

From (24) and (25), it follows $\frac{n-2}{2}\frac{n^{2}}{n-1}\left\|
H\right\| ^{2}\leq \frac{n^{2}}{2}\frac{2n-3}{2n+3}\left\| H\right\| ^{2}$,
whence\newline
$\frac{n^{2}}{2}\left\| H\right\| ^{2}(\frac{2n-3}{2n+3}-\frac{n-2}{n-1}%
)\geq 0.$

As $\frac{2n-3}{2n+3}-\frac{n-2}{n-1}<0,$ $\forall $ $n\geq 3,$ we infer
that $\left\| H\right\| =0$, so $M$ is a minimal submanifold in $\widetilde{M%
}(c).$

\ 

{\bf Acknowledgment. }I would like to thank to Professor C. Udri\c{s}te, who
has always been generous with his time and advice.

\end{document}